\documentclass[12pt]{amsart}

\usepackage{amssymb,amsmath,amsfonts}  
\usepackage{verbatim,a4wide}
\usepackage[cp1250]{inputenc}
\usepackage[OT4]{fontenc}
\usepackage[english]{babel}

\usepackage{bm}
\usepackage{graphicx}
\usepackage{arydshln}
\usepackage{amsfonts}
\usepackage{color}

\numberwithin{equation}{section}
\newtheorem{thm}{Theorem}[section]

\newtheorem{alg}[thm]{Algorithm}
\newtheorem{prob}[thm]{Problem}

\newenvironment{exmp}{\refstepcounter{thm}\ \\[2mm]%
               \noindent{\bf Example~\thethm\ }}{ \qed\\\rm}
\newenvironment{pf}{{\it Proof.\/}\rm }{\qed\\[-1ex]}

\renewcommand{\arraystretch}{1}
\newcommand{\R}{\mathbb R}

\newcommand{\ip}[2]{\mbox{$\left\langle{#1},{#2}\right\rangle$}}
\newcommand{\dt}{\mbox{$\der t$}}
\newcommand{\dx}{\mbox{$\der x$}}
\newcommand{\der}{\mbox{${\rm\,d}$}}
\newcommand{\const}{\textrm{\,const}}

\newcommand{\rbinom}[2]{\mbox{$\displaystyle\binom{#1}{#2}^{\!\!-1}\!\!$}}

\begin{document}

	\title[Weighted polynomial approximation of rational B\'ezier curves]
	{Weighted polynomial approximation\\ of rational B\'ezier curves${}^\ast$}
	\author{Stanis{\l}aw Lewanowicz${}^\dagger$, Pawe{\l} Wo\'zny, Pawe{\l} Keller}
\address{Institute of Computer Science, University of Wroc{\l}aw, Poland}
\email{\{Stanislaw.Lewanowicz,Pawel.Wozny\}@ii.uni.wroc.pl}

\address{Faculty of Mathematics and Information Science,  Warsaw University of Technology, Poland}
\email{Pawel.Keller@mini.pw.edu.pl}
\date{\today}

\thispagestyle{empty}
\maketitle
\footnotetext[2]{Corresponding author. Fax + 48 71 3757801.}
\footnotetext[1]{This paper is an extended version of our paper \cite{LWK12}, 
	in which the simplest form  of the distance between curves is used.}
\noindent \textbf{Abstract}.
 We present an efficient method to solve the problem of the
constrained least squares  approximation
 of the rational B\'ezier curve  by the B\'ezier curve. 
The presented algorithm uses the dual constrained Bernstein basis polynomials, 
associated with the Jacobi scalar product, and exploits their recursive properties.
Examples are given,  
showing the effectiveness of the algorithm.\\

\noindent \textit{AMS classification}: {Primary 41A10. Secondary 65D17, 33D45}

\noindent \textit{Keywords}: Rational B\'ezier curve; Polynomial approximation;
                  Constrained dual  Bernstein basis.                    



\section{Introduction}
                                                        \label{SS:Intr}
In CAGD, it is frequently important to approximate a rational B\'ezier curve by
a polynomial one. In the last two decades, many approaches to this problem have 
been proposed \cite{SK91,WSC97,LW00,Flo06,LP98,Lu11,HSL08}. The large spectrum
of methods contains approximation by hybrid curves \cite{SK91,WSC97,LW00},
Hermite interpolation \cite{Flo06,LW00}, progressive iteration approximation \cite{Lu11},
least squares approximation \cite{LP98} and approximation by B\'ezier curves
with the control points obtained by successive degree elevation of the rational
B\'ezier curve \cite{HSL08}.

In this paper, we consider   the following approximation problem.
\begin{prob}\label{P:main}
Let be given a rational B\'ezier curve  $R_n$ of  degree $n$, with control points $r_i\in\R^d$ and positive weights $\omega_i\in\R$,
\begin{equation}
	\label{E:Rn-intr}
	R_n(t)=\dfrac{Q_n(t)}{\omega(t)}=\frac{ \sum_{i=0}^n \omega_ir_i\,B^n_i(t)}
	                                                               {  \sum_{i=0}^n \omega_i\,B^n_i(t)}
\qquad (0\le t\le 1),
\end{equation}
where 
\[
   B^n_{i}(t)=\binom ni t^i(1-t)^{n-i} \qquad (0\le i\le n)
\]
are {Bernstein basis polynomials}. 
Find a degree $m$  B\'ezier curve
\begin{equation}
	\label{E:Pm-intr}
	P_m(t)=\sum_{i=0}^m p_i\,B^m_i(t)\qquad (0\le t\le 1)
\end{equation}
such that the value of the error
\begin{equation}
          \label{E:sqerr}
	\int_{0}^{1}(1-t)^\alpha t^\beta\,\|R_n(t)-P_m(t)\|^2\dt\qquad (\alpha,\;\beta>-1)
\end{equation}
is minimized in the space $\Pi^d_m$ of parametric polynomials in $\R^d$ of degree at most $m$
(for simplicity, we write $\Pi_m:=\Pi^1_m$)
under the additional conditions that
\begin{equation}\label{E:constraints}
	\left\{\begin{array}{l}
		R_n^{(i)}(0)=P_m^{(i)}(0) \qquad (i=0,1,\ldots,k-1),\\[2ex]   
		R_n^{(j)}(1)=P_m^{(j)}(1) \qquad (j=0,1,\ldots,l-1),  
	\end{array}\right.
\end{equation}
where $k+l\leq m$. Here $\|\cdot\|$ is the Euclidean vector norm.
\end{prob}
Note that in the case $\alpha=\beta=0$,  the above problem as well the method proposed in this paper reduce to the form given in \cite{LWK12}.

The values of $n$ and $m$ are not related. However, if $(m+1)d<(n+1)(d+1)$ then the number 
of parameters of the approximating polynomial curve is smaller than the total number 
of parameters of the rational curve.

In \cite{LP98},  an \textit{approximate} solution to the above problem with $m>n$ 
is obtained by solving a linear least squares problem. Paper \cite{CW10} deals with
more general problem of the constrained degree reduction of rational B\'ezier curves; 
one of the auxiliary problems discussed there contains the above problem as a particular case.
We present a method which is based on the idea of using constrained dual Bernstein polynomial basis  
\cite{WL09} to compute the control points $p_i$. Our algorithm is  efficient thanks
to  using  fast schemes of  $1^o$ evaluation  the B\'ezier form coefficients of the dual 
polynomials \cite{LW11b} and  $2^o$~numerical computation of the collection of integrals
\[
	\int_{0}^{1}\frac{(1-t)^\alpha t^\beta B^{n+m}_h(t)}{\omega(t)}\dt \qquad (k\le h\le n+m-l).
\]
The cost  is significantly lower than in the case of the special variant of the 
method of \cite{CW10}, which needs inverting a $(n+m-k-l+1)\times(n+m-k-l+1)$ matrix.  

Let us mention that the problem stated above can be also considered for other 
norms. However, even the simpler problem of constrained degree reduction
of B\'ezier curves in $L_\infty$-norm 
requires higher computational complexity (see, e.g., \cite{Ahn03}). 
The most appropriate metric for curves in geometric terms would be the Hausdorff distance, but
the computation of such distance of the nonlinear curves is not so easy.
Hence, $L_2$-norm seems to be a good choice as we can construct solution 
in explicit form using the Bernstein and dual Bernstein bases in a natural 
and convenient way.  

The outline of this paper is as follows. Section~\ref{S:constrdual} contains 
basic facts on the constrained dual Bernstein polynomials.
 Section~\ref{S:appr} brings a complete solution to Problem~\ref{P:main};
 for implementation details, see Section~\ref{S:impl}.
 In Section~\ref{S:Exmp}, the proposed method is applied to some examples
 and compared with two other algorithms.

\section{Constrained dual Bernstein polynomials} \label{S:constrdual}

Let $\Pi_m^{(k,l)}, $where $k+l\le m$, be the space of all polynomials of degree $\le m$, 
whose derivatives of order  $\le k-1$ at $t=0$, as well as derivatives of order
$\le l-1$ at $t=1$, vanish:
\[
\Pi_m^{(k,l)}:=
       \left\{P\in\Pi_m\::\:
       P^{(i)}(0)=0\quad (0\le i\le k-1)
       \; \mbox{and}\;        
       P^{(j)}(1)=0\quad (0\le j\le l-1)\right\}.
\]
 Obviously,  $\mbox{dim}\;\Pi_m^{(k,l)}=m-k-l+1$, and the Bernstein polynomials
\(
\left\{B^m_k,B^m_{k+1},\ldots,B^m_{m-l}\right\}
\)
form a basis of this space.
There is a unique \textit{dual constrained Bernstein basis of degree} $m$,
\[
D^{(m,k,l)}_k(x;\alpha,\beta),\,D^{(m,k,l)}_{k+1}(x;\alpha,\beta),\ldots,
D^{(m,k,l)}_{m-l}(x;\alpha,\beta),
\]
satisfying
\[
	\langle D^{(m,k,l)}_i,\,B^m_j \rangle=\delta_{ij}\qquad (i,j=k,k+1,\ldots,m-l),
\]
where $\delta_{ij}$ is 1 if $i=j$ and 0 otherwise, and the inner product $\langle\cdot,\cdot\rangle$ is given by
\[
	\langle f,\,g \rangle:=\int_{0}^{1}(1-t)^\alpha t^\beta f(t)g(t)\dt\qquad
	(\alpha,\,\beta>-1).
\]

Properties of  the  polynomials $D^{(m,k,l)}_i$ are studied in  \cite{LW11a,LW11b} and
\cite{WL09}. 
We need the following result.
\begin{thm}[\cite{LW11b}] \label{T:constrDinB}
The constrained dual basis polynomials have the  B\'ezier-Bernstein representation
\begin{equation}
	\label{E:constrDinB}
	 D^{(m,k,l)}_i(x;\alpha,\beta)=\sum_{j=k}^{m-l}c_{ij}(m,k,l,\alpha,\beta)\,  
	                   B^{m}_j(x),
\end{equation}
where 
the coefficients $c_{ij}\equiv c_{ij}(m,k,l,\alpha,\beta)$ 
satisfy the recurrence relation
\begin{align}\label{E:c-rec}
c_{i+1,j}=&\frac1{A(i)}\,\left\{(i-j)(2i+2j-2m-\alpha+\beta)\,c_{ij}\right. \\[1.5ex]
&\hphantom{\frac1{A(i)}\,\{}
\left.+B(j)\,c_{i,j-1} +A(j)\,c_{i,j+1}-B(i)\,c_{i-1,j}\right\}\nonumber \\[1.5ex]
&\hphantom{\frac1{A(i)}\,\{+B(j)\,c_{i,j-1} +A(j)\,c_{i,j+1}}
(k\le i\le m-l-1,\quad  k\le j\le m-l) \nonumber
\end{align}	
with
\[ 
\begin{array}{l}
A(u):=(u-m)(u-k+1)(u+k+\beta+1)/(u+1),\\[2ex]	
B(u):=u(u-m-l-\alpha-1)(u-m+l-1)/(u-m-1).
\end{array}
\]  
We adopt the convention that
$c_{ij}:=0$ if $i<k$, or $i>m-l$, or $j<k$, or $j>m-l$.
The starting values are
\begin{align}	\label{E:Crec-start}
c_{kj}
 =&\rbinom{m}{k}\,
 \frac{(-1)^{k}(\sigma+2k+2l+1)_{m-k-l}(k+\beta+2)_{m-l}}{(m-k-l)!B(\alpha+2l+1,\,\beta+2k+1)} \\[2ex]
 &\times \binom{m-k-l}{j-k}\rbinom{m}{j}\,\,\frac{(-1)^j}{(\alpha+2l+1)_{m-l-j}(k+\beta+2)_j},  \nonumber	
\end{align}
where $j=k,k+1,\ldots, m-l$, and we use the notation  
\[
	(a)_k:=\prod_{j=0}^{k-1}(a+j)\quad (k\ge0),\qquad
\sigma:=\alpha+\beta+1,\qquad B(\lambda,\,\mu)=\frac{\Gamma(\lambda)\Gamma(\mu)}{\Gamma(\lambda+\mu)}.
\]
\end{thm}
Observe that the quantities $c_{ij}$ can be put in a square table (see~Table~\ref{Tab:C}).
\begin{table}[htb]
\caption{The $c$-table \label{Tab:C}}
\normalsize
\[
\begin{array}{cccccc}
    &0&0&\ldots&0& \\[2ex]
    0& c_{kk}&c_{k,k+1}&\ldots&c_{k,m-l}&0     \\[2ex]
    0& c_{k+1,k}&c_{k+1,k+1}&\ldots&c_{k+1,m-l}&0  \\[2ex]
    \multicolumn{6}{c}{\dotfill}\\[2ex]
    0& c_{m-l,k}&c_{m-l,k+1}&\ldots&c_{m-l,m-l}&0  \\[2ex]
     &0&0&\ldots&0& 
\end{array}
\]	
\end{table}
Now,   the $c$-table can be completed very easily in the following way.
\begin{alg}[Computing the coefficients $ c_{ij}(m,k,l,\alpha,\beta)$  \cite{LW11b}]
        \label{A:ctab}
        \ \begin{description}
        \item[I] Compute recursively quantities $c_{kk},\;c_{k,k+1},\ldots,\,c_{k,m-l},$
filling the first row of the $c$-table, by the formulas
\begin{align}
	\label{E:ckml}
	  c_{k,m-l}:=&\rbinom{m}{k}\,\rbinom{m}{l}\;
	  \frac{(-1)^{m-k-l}(\sigma+2k+2l+1)_{m-k-l}}
	              {B(\alpha+2l+1,\,\beta+2k+1)\,(m-k-l)!},\\[2ex]
	\label{E:Ckj}
	  c_{kj}:=&\frac{(j-m)(j-k+1)(j+\beta+k+2)}{(j+1)(j-m+l)(j-\alpha-l-m)}\,c_{k,j+1} \\[1ex]
	  &\hphantom{\frac{(j-n)(j-k+1)(j+\beta+k+2)}{(j+1)(j-n+l)(j-\alpha-l-n)}} 
	  (j=m-l-1,m-l-2,\ldots,k). \nonumber
\end{align}
\item[II] For $i=k,k+1,\ldots,m-l-1$ and $j=k,k+1,\ldots,m-l$, compute
      			$c_{i+1,j}$, using the recurrence   \eqref{E:c-rec}.
      		\end{description}
\end{alg}
\section{Constrained polynomial approximation of rational B\'ezier curve}
						\label{S:appr}
Clearly, the B\'ezier curve being the solution of Problem~\ref{P:main} can be obtained 
in a componentwise way. Hence, it is sufficient to give the details of our method of solving
this problem in case where $R_n,\;P_m\in\R^1$.

Given the rational function  
\begin{equation}
	\label{E:Rn}
	R_n(t)=\dfrac{Q_n(t)}{\omega(t)}
	=\frac{ \sum_{i=0}^n \omega_ir_i\,B^n_i(t)}
	          {  \sum_{i=0}^n \omega_i\,B^n_i(t)}
\qquad (0\le t\le 1),
\end{equation}
with  $r_i\in\R$
and  $\omega_i\in\R_+$,
 we look for a degree $m$  polynomial
\begin{equation}
	\label{E:Pm}
	P_m(t)=\sum_{i=0}^m p_i\,B^m_i(t)\qquad (0\le t\le 1)
\end{equation}
which gives the minimum value of
\begin{equation}\label{E:err}
	\|R_n(t)-P_m(t)\|^2_{L_2}=\langle R_n-P_m,R_n-P_m\rangle
\end{equation}
with the constraints
\begin{equation}\label{E:constr}
	\begin{array}{l}
		R_n^{(i)}(0)=P_m^{(i)}(0) \qquad (i=0,1,\ldots,k-1),\\[2ex]   
		R_n^{(j)}(1)=P_m^{(j)}(1) \qquad (j=0,1,\ldots,l-1),  
	\end{array}
\end{equation}
where $k+l\leq m$.
\begin{thm}
	\label{T:main}
Given the coefficients $r_0,r_1,\ldots,r_n$ and weights $\omega_0,\omega_1,\ldots,\omega_n$ of 
the rational function \eqref{E:Rn}, the coefficients $p_0,p_1, \ldots,p_m$ of the polynomial
\eqref{E:Pm} minimising the error \eqref{E:err} with constraints \eqref{E:constr}
are given by
\begin{align}
	\label{E:pi-begin}
	  p_i=&\frac{(m-i)!}{m!}\varrho_{i0}-\sum_{j=0}^{i-1}(-1)^{i+j}\binom{i}{j}p_j
	  \qquad (i=0,1,\ldots,k-1); \\
	\label{E:pi-end}
	  p_{m-i}=&(-1)^i\frac{(m-i)!}{m!}\varrho_{i1}
	          -\sum_{j=1}^{i}(-1)^{j}\binom{i}{j}p_{m-i+j}
		   \qquad (i=0,1,\ldots,l-1);\\
	 \label{E:pi-mid}
	  p_{i}=&\sum_{h=0}^{n}\binom{n}{h}\omega_h\,r_h\,
	  \sum_{j=k}^{m-l}
	  \binom{m}{j}\rbinom{n+m}{j+h}\,c_{ij}(m,k,l,\alpha,\beta)
	  I_{j+h} \\
	  &-
	  \left(\sum_{j=0}^{k-1}+\sum_{j=m-l+1}^{m}\right)p_j\,K_{ij}
	  \qquad 
	  (i=k,k+1,\ldots,m-l), \nonumber	
\end{align}
where  $c_{ij}(m,k,l,\alpha,\beta)$ are introduced in \eqref{E:constrDinB}, while 
$\varrho_{i0}$ and $\varrho_{i1}$  are defined recursively for $i=0,1,\ldots$ by
\begin{align}
	\label{E:rho0}
	  \varrho_{i0}:=&\frac{n!}{\omega_0}\left[\frac{1}{(n-i)!}
	     \Delta^i\left(\omega_0r_0\right)-\sum_{j=0}^{i-1}\binom{i}{j}
	     \frac{\Delta^{i-j}\omega_0}{(n-i+j)!}\,\varrho_{j0}\right];\\
	\label{E:rho1}
	  \varrho_{i1}:= &\frac{n!}{\omega_n}\left[\frac{1}{(n-i)!}
	     \Delta^i\left(\omega_{n-i}r_{n-i}\right)-\sum_{j=0}^{i-1}\binom{i}{j}
	     \frac{\Delta^{i-j}\omega_{n-i+j}}{(n-i+j)!}\,\varrho_{j1}\right],
\end{align}
and
\begin{align}
	\label{E:I}
	I_h:=&\int_{0}^{1}\frac{(1-t)^\alpha t^\beta B^{n+m}_h(t)}{\omega(t)}\,\dt,\\[1ex]	
	\label{E:K}
	K_{ij}	
	=&\binom{m}{j}\rbinom{m}{i}\,\frac{(-1)^{i-k}(k-j)_{m-k-l+1}}{(i-j)(i-k)!(m-l-i)!}
       \\[1ex]
	&\times
	\frac{(\alpha+l+1)_{m-j}(\beta+k+1)_{j}}{(\alpha+l+1)_{m-i}(\beta+k+1)_{i}}.\nonumber	  
\end{align}
Here we use the standard notation $\Delta^0c_h=c_h$, $\Delta^j
c_h=\Delta^{j-1}c_{h+1}-\Delta^{j-1}c_{h}$ ($j=1,2,\ldots$). 
\end{thm}

\begin{pf} Recall that for arbitrary polynomial of degree $N$,  
\[
	U_N(t)=\sum_{i=0}^{N}u_i\,B^N_i(t),
\]
the well-known formulas hold (see, e.g., \cite[p. 49]{Far96})
\begin{eqnarray*}
U^{(j)}_N(0)&=&\frac{N!}{(N-j)!}\Delta^ju_0=
\frac{N!}{(N-j)!}\sum_{h=0}^{j}(-1)^{j+h}\binom{j}{h}u_{h},\\
U^{(j)}_N(1)&=&\frac{N!}{(N-j)!}\Delta^j u_{N-j}=
\frac{N!}{(N-j)!}\sum_{h=0}^{j}(-1)^{j+h}\binom{j}{h}u_{N-j+h}. 
\end{eqnarray*}
Using them in
\[  
 R^{(i)}_n(h)=\frac{1}{\omega(h)}\left[Q^{(i)}_n(h)-\sum_{j=0}^{i-1}
 \binom{i}{j}R^{(j)}_n(h)\omega^{(i-j)}(h)\right]\quad (h=0,1)	
\]  
gives equations \eqref{E:rho0}, \eqref{E:rho1}, where we denoted
\[
	\varrho_{ih}:=R^{(i)}_n(h).
\]
Using the above equations 
in \eqref{E:constr}, we obtain the forms 
\eqref{E:pi-begin} and
\eqref{E:pi-end} for the coefficients $p_0,\,p_1,\ldots,p_{k-1}$ and
$p_{m-l+1},\ldots,p_{m-1},p_m$, respectively.

The remaining coefficients $p_i$ are to be determined so that
\[
\| R_n-P_m\|^2_{L_2}=\| W-\sum_{i=k}^{m-l}p_iB^m_i\|^2_{L_2}	
\]
has the least value, where
\[
W:=R_n-\left(\sum_{j=0}^{k-1}+\sum_{j=m-l+1}^{m}\right)p_jB^m_j.
\]

Remembering that  $B^m_i$ and $D^{(m,k,l)}_i$ ($k\le i\le m-l$) are dual
bases in the space $\Pi_m^{(k,l)}$, we obtain the formula
\begin{align*}
p_i=&\ip{W}{ D^{(m,k,l)}_i}
\\=&
\sum_{j=k}^{m-l}c_{ij}(m,k,l,\alpha,\beta)\,\langle R_n,\,B^m_j\rangle 
    -\left(\sum_{j=0}^{k-1}+\sum_{j=m-l+1}^{m}\right)p_j
  \ip{B^m_j}{D^{(m,k,l)}_i}.
\end{align*} 
Observe that
\begin{align*}
	\langle R_n,\,B^m_j\rangle=&\sum_{h=0}^{n}\omega_h\,r_h
	\langle {B^n_h}/{\omega},\,B^m_j\rangle\\
	=&\sum_{h=0}^{n}\omega_h\,r_h
	\binom{n}{h}\binom{m}{j}\rbinom{m+n}{h+j}
	\langle 1/{\omega},\,B^{m+n}_{h+j}\rangle\\
	=&\sum_{h=0}^{n}\omega_h\,r_h
	\binom{n}{h}\binom{m}{j}\rbinom{m+n}{h+j}
	I_{h+j},
\end{align*}
where we use notation \eqref{E:I}.

Using results of \cite{WL09} and \cite{LW11a}, we deduce that 
\begin{align*}
	K_{ij}:=&\ip {B^m_j}{D^{(m,k,l)}_i}\\
=&\binom{m}{j}\rbinom{m}{i}\,
   \frac{(\alpha+l+1)_{m-j}(\beta+k+1)_{j}}{(i-k)!(m-l-i)!(\alpha+l+1)_l(\beta+k+1)_k}\\[1ex]
	&\times d_{i-k}(j-k;\beta+2k,\alpha+2l,m-k-l),
\end{align*}
where for $h=0,1,\ldots, N$ we define
\begin{align*}
	&d_h(x;a,b,N)=  \lim_{n\to N}\,\frac{(-1)^h\,N!}{n!(a+1)_h(b+1)_{n-h}} 
 \,\sum_{v=0}^{h}\frac{(-h)_v\,(v+1-N)_{n-v}(b+1)_{n-v}}
                      {(n+a+b+2)_{N-v}}  \\ 
  &\hphantom{d_h(x;a,b,N)= \lim_{n\to N}\,\times \sum_{v=0}^{h}\frac{(-h)_v\,(b+1)_{n-v+2}}{a}}
  \times 
  \,Q_{n-v}(N-x;b,a+v+1,N-v-1)
\end{align*}
with 
\[
	Q_p(t;\mu,\nu,M):=
	\sum_{u=0}^{p}\frac{(-p)_u(p+\mu+\nu+1)_u(-t)_u}{u!(\mu+1)_u(-M)_u}
\]
being Hahn polynomials (see, e.g., \cite[\S9.5]{KLS10}).
It is easy to see that
\[
	d_h(x;a,b,N)=\frac{(-x)_{N+1}(-N-b)_h}{(h-x)(a+1)_h(b+1)_N}\qquad (0\le h\le N).
\]
 Hence, some algebra gives
\begin{align*}		
	K_{ij}	
	=&\binom{m}{j}\rbinom{m}{i}\,\frac{(-1)^{i-k}(k-j)_{m-k-l+1}}{(i-j)(i-k)!(m-l-i)!}
      	\nonumber \\[1ex]
	&\times
	\frac{(\alpha+l+1)_{m-j}(\beta+k+1)_{j}}{(\alpha+l+1)_{m-i}(\beta+k+1)_{i}},
\end{align*}
which is formula \eqref{E:K}.

Now, \eqref{E:pi-mid} readily follows.
\end{pf}

\section{Implementation of the method}
				\label{S:impl}

\subsection{Computing integrals \eqref{E:I}}
				\label{SS:int}
Integrals \eqref{E:I} involving rational function cannot be evaluated exactly.
However, we show that they can be computed numerically up to high precision
using the method described in \cite{Kel07}.

Observe that formula \eqref{E:I} can be written as
\begin{equation}
	\label{E:Ivar}
I_h=2^{-\sigma-N}\binom{N}{h}J(\alpha+N-l-h,\beta-k+h;\vartheta)
\qquad (k\le h\le N-l;\;N:=n+m),	
\end{equation}
where 
\begin{equation}
	\label{E:J}
J(a,b;\vartheta):=\int_{-1}^{1}(1-x)^{a}(1+x)^{b}\vartheta(x)\dx,	
\end{equation}
and
\[
\vartheta(x):=\frac{(1-x)^l(1+x)^k}{\omega\left(\tfrac12(1+x)\right)}.
\]
Notice that
by assumption on the positivity of the weights $\omega_i$, the polynomial $\omega(t)$ has no roots in the interval $[0,\,1]$, hence the function 
$\vartheta(x)$ is analytic in a planar region containing the interval $[-1,\,1]$. This implies that 
the function $\vartheta$  can be well approximated  by a sum of Chebyshev  polynomials $T_j(x)$ of the first
kind:
\begin{equation}
	\label{E:SM}
 S_M(x)=\frac12 \gamma_0\,T_0(x)+\sum_{j=1}^{M} \gamma_j\,T_j(x)\qquad (-1\le x\le 1).  
\end{equation}
See, e.g., \cite{Gen72}. 
Clearly,
\[
	J(a,b;\vartheta)\approx J(a,b;S_M).
\]
The coefficients $\gamma_j$ are determined so that $S_M$ interpolates $\vartheta$
at the abscissae $\xi_i:=\cos(i\pi/M)$ ($0\le i\le M$), hence
\begin{equation}
	\label{E:gam}
 \gamma_j:=\frac{2-\delta_{jM}}{M}	\sum_{i=0}^{M}{}^{\!''}\vartheta\left(\cos\frac{i\pi}{M}\right)\cos\frac{ij\pi}{M}
 \qquad(0\le j\le M).
\end{equation}
(The double prime on the sum  means that the first and the last terms 
are to be halved.)
The right-hand side of \eqref{E:gam} is known to be efficiently computed by means of the FFT for real data
(see \cite{Gen72}, or \cite[\S5.1]{DB08}; the authors  recall that the FFT is not
only fast, but also resistant to roundoff errors). In practical implementation, 
the coefficients 
$\gamma_j$ $(0\leq j\leq M)$ are computed repeatedly
for doubled values of $M$ $(M=32,64,\dots)$ until
\begin{equation}\label{E:Mcrit}
  \sum\limits_{i=M-3}^{M} |\gamma_i| < \varepsilon,
\end{equation}
where  $\varepsilon$ is a prescribed tolerance.
 For better efficiency, in every step one should reuse all the previously
computed values of the function $\vartheta$.

Now, we have the following result.
\begin{alg}[Numerical computation of the integral $J(a,b;S_M)$]
			\label{A:polint}
Given $a,\,b>-1$, let  $r:=b-a$ and  $s:=a+b+1$.
Define the sequence $d_i$ ($0\le i\le M+1$) 
recusively by
\begin{equation}\label{E:di}
	  \begin{array}{l}
	  d_{M+1}=d_M:=0,  \\[1ex]
	  d_{i-1}:=\dfrac{2r\,d_i+(i-s)\,d_{i+1}-2\,\gamma_i}{i+s}
	  \qquad (i=M,M-1,\ldots,1).
	  \end{array}
\end{equation}
Then we have
\begin{equation}
	\label{E:polint}
J(a,b;S_M) = 2^{s-1}B(a+1,b+1)\,(\gamma_0-r\,d_0+s\,d_1).	
\end{equation}
\end{alg}
\begin{pf}
Function $\rho(x)=(1-x)^a(1+x)^b$ is a solution of the differential equation
\[
	(1-x^2)\rho'(x)-[r-x(s-1)]\rho(x)=0,
\]
where $r:=b-a$, $s:=a+b+1$.
By \cite[Thm 2.1]{Kel07}, we have
\[
	\int (1-x)^a(1+x)^bS_M(x)\dx=(1-x)^{a+1}(1+x)^{b+1}\,T(x)+\const
	\qquad (-1\le x\le 1)
\]
for any continuous solution $T$ of the differential equation
\begin{equation}
	\label{E:diffeq}
(1-x^2)T'(x)+[r-x(s+1)]T(x)=S_M(x)\qquad (-1\le x\le 1).	
\end{equation}
Further, using the approach of \cite{Lew92} yields the recurrence relation
\[
	(-i-s)\,d_{i-1}+2r\,d_i+(i-s)\,d_{i+1}=2\,\gamma_i\qquad (i=0,1,\ldots)
\]
with $d_{-1}:=d_1$
	and $\gamma_j=0$ for $j>M$, satisfied by the Chebyshev coefficients $d_i$ of the solution $T$ of equation \eqref{E:diffeq}.

However, equation \eqref{E:diffeq} has, in general, no continuous solution on $[-1,\,1]$.
On the other hand, by \cite[Thm 2.3]{Kel07} there always exists a \textit{polynomial} solution
of the modified equation
\begin{equation}\label{E:modeq}
(1-x^2)T'(x)+[r-x(s+1)]T(x)=S_M(x)-\tfrac12\delta\qquad (-1\le x\le 1),
\end{equation}
where
\(
	\delta:=\gamma_0-r\,d_0+s\,d_1.
\)
Obviously,
\[
	J\left(a,b;S_M-\tfrac12\delta\right)=\int_{-1}^1 (1-x)^a(1+x)^b[S_M(x)-\tfrac12\delta]\dx=0,
\]
hence
\[
	J\left(a,b;S_M\right)=J\left(a,b;\tfrac12\delta\right)=\tfrac12\delta\,2^sB(a+1,b+1).
\]		
\end{pf}

\subsection{Main algorithm}
				\label{SS:Alg}
The presented method is summarized in the following algorithm. 
\begin{alg}[Constrained polynomial approximation of the rational B\'ezier curve]	
				\label{A:main}
Given the coefficients $r_0,r_1,\ldots,r_n$ and $\omega_0,\omega_1,\ldots,\omega_n$ of 
the rational function 
\[
R_n(t)=\dfrac{Q_n(t)}{\omega(t)}
	=\frac{\sum_{i=0}^n \omega_ir_i\,B^n_i(t)}{ \sum_{i=0}^n \omega_i\,B^n_i(t)}
		\qquad (0\le t\le 1),
\]
the coefficients $p_0,p_1, \ldots,p_m$ of the polynomial
\[
		P_m(t)=\sum_{i=0}^m p_i\,B^m_i(t)\qquad (0\le t\le 1)
\]
minimising the error \eqref{E:err} with constraints \eqref{E:constr}
are computed in the following way.\\[-2ex]
\begin{description}
\itemsep0pt
\item[{\sc Step 1}] Compute $p_0,\,p_1,\ldots,p_{k-1}$ by \eqref{E:pi-begin}.

\item[{\sc Step 2}] Compute $p_m,\,p_{m-1},\ldots,p_{m-l+1}$ by \eqref{E:pi-end}.

\item[{\sc Step 3}] Compute  $c_{ij}(m,k,l,\alpha,\beta)$
                    for $i,j=k,\ldots,m-l$ by Algorithm~\ref{A:ctab}.
		    
\item[{\sc Step 4}] Given $\varepsilon>0$,   compute the coefficients
			$\gamma_j$
                      of the polynomial $S_M$ (cf. \eqref{E:SM})
                      by using FFT, with $M$ determined so 
                       that \eqref{E:Mcrit} holds.
		      
\item[{\sc Step 5}] For $h=k,k+1,\ldots, N-l$, where $N:=m+n$,
                 compute
		\[
		J^M_h:=J(\alpha+N-l-h,\,\beta-k+h;\,S_M)
		\]
		 by Algorithm~\ref{A:polint}.

\item[{\sc Step 6}] For $i=k,k+1, \ldots,m-l$, compute $p_i$ 
             by \eqref{E:pi-mid} with $I_h$ replaced by
             \[
             I^M_h:=2^{-\sigma-N}\binom{N}{h}\,J^M_h\qquad (k\le h\le N-l).
             \]

\end{description}
\end{alg}

\section{Examples}
                     \label{S:Exmp}
In this section, we present several examples of approximation of rational B\'ezier
surfaces by B\'ezier curves with constraints, which we have described in Section~\ref{S:appr}.
Computations were carried out on a computer with
Intel Core i5 3.33GHz processor and 4GB of RAM.
We used 18-digit arithmetic and set 
$\varepsilon:=10^{-16}$ in  Step 4 of the Algorithm~\ref{A:main} to ensure that the integrals \eqref{E:I} are computed within the accuracy close to the representation error.

In the examples below, we use the notation
\[
	E(t):=\|R_n(t)-P_m(t)\|
\]
for the error function, and 
\[ 
	e_\infty:=\max_{t\in[0,\,1]}E(t), \qquad
	 e_2(\alpha,\beta):= \left(\int_{0}^{1}(1-t)^\alpha t^\beta E^2(t)\,\dt\right)^{1/2}
\] 
for the maximum and least-squares approximation error, respectively.

\begin{exmp}[\textit{Starling's sketch}]
\label{Ex:1}
We consider the curve shown in Figure~\ref{fig:Fig1}  (solid and red),
obtained by joining  two rational B\'ezier curves of degree eight each; we say that
the composite curve has degree (8,8).  
The first  curve is defined by the control points 
 $(23, 57)$,  $(-13, 43)$,  $(29, 58)$,  $(44, 48)$,  $(30, 42)$, 
             $(13, 44)$,  $(-2, 77)$,  $(42, 83)$,  $(80, 1)$,
  and the associated weights
  $1,\, 4,\, 3,\, 1,\, 5,\, 4,\, 7,\, 6,\, 1$, while
the second one by $(80, 1)$,  $(14, 4)$, $(3, 54)$, $(42, 54)$,  $(51, 42)$,  
             $(36, 49)$,  $(66, 12)$,  $(36, 2)$,  $(47, 3)$,
and $1,\, 1,\,  4,\,  4,\,  2,\, 3,\, 3,\, 7,\, 8$. 
The polynomial composite curve approximation of the above curve with $\alpha=\beta=\frac12$,
constraints of $C^1$- and $C^2$-continuity, 
  without and with subdivision, is shown in Figure~\ref{fig:Fig1}.
  	\begin{figure}[htbp]
	\begin{center}
        \hspace*{-1cm} 	
	\includegraphics[scale=0.45,angle=0]{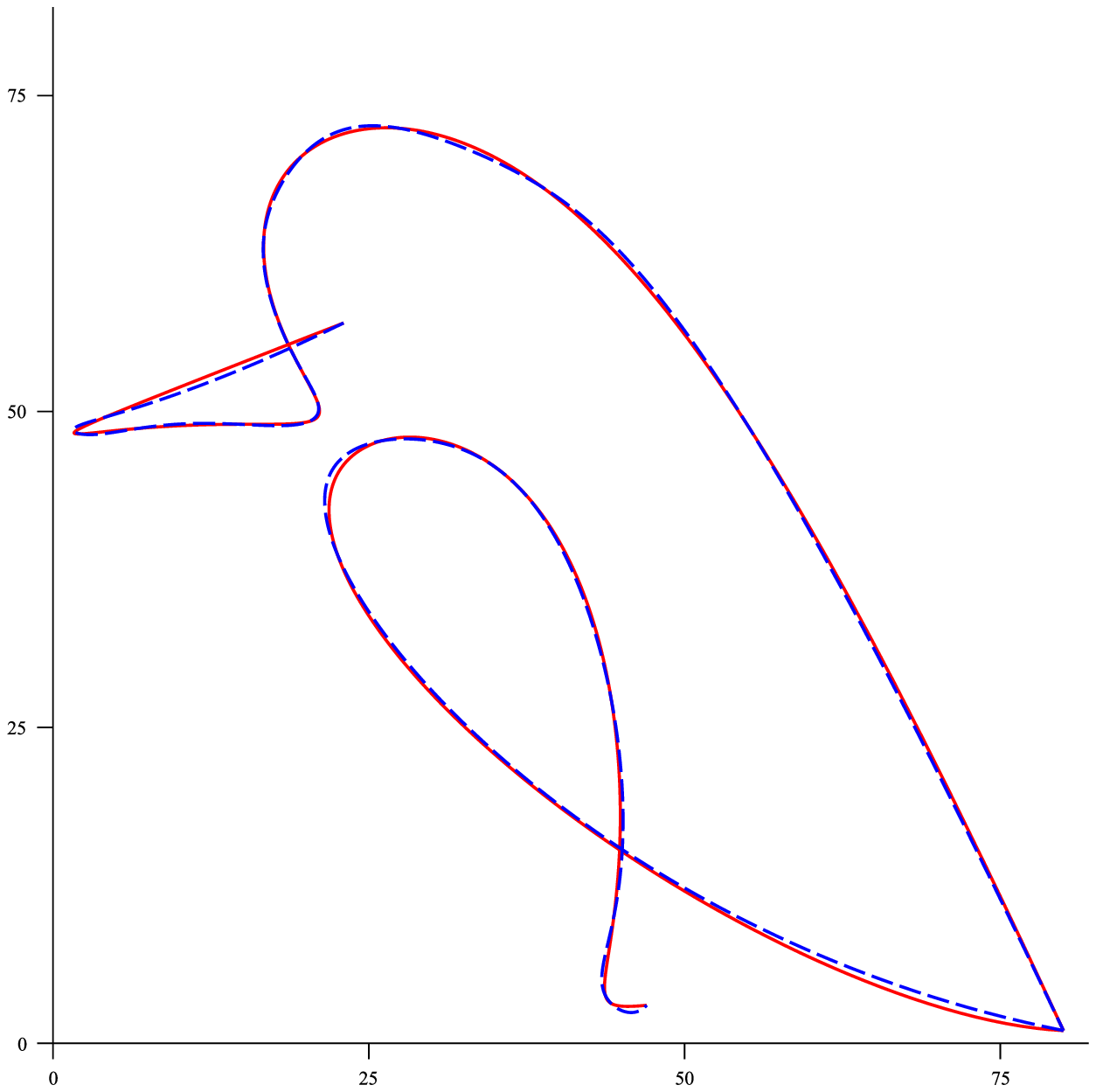}  
	\hspace*{-0.75cm} 	
	\includegraphics[scale=0.29,angle=0]{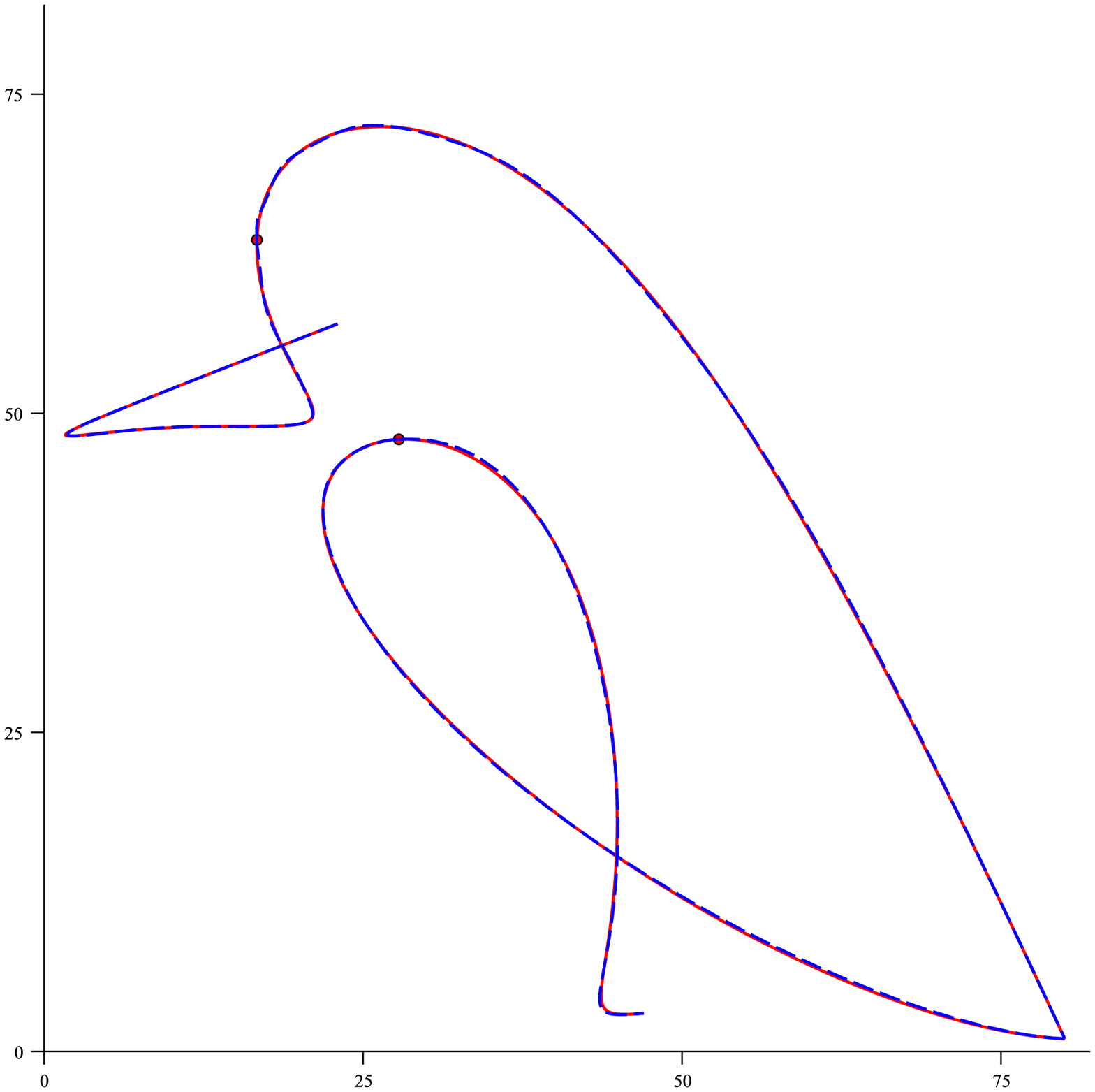} 
        \caption{\textit{Left}: B\'ezier composite curve  approximation (dashed and blue) of degree (13,8) 
        to the rational composite B\'ezier
	curve (solid and red) of degree (8,8) with the end-point interpolation ($k=l=1$). 
	The errors are $e_\infty=\{3.152,\,2.814\}$ and
	$e_2(\frac12,\frac12)=\{0.166,\,0.284\}$, respectively.
	\textit{Right}: B\'ezier composite curve  approximation (dashed and blue) of degree (12,11,7,6)
        to the same rational composite B\'ezier
	curve (solid and red) of degree (8,8), with one subdivision of each of its two parts, 
	and the end-point derivative interpolation ($k=l=2$). 
	The errors are $e_\infty=\{0.559,\,0.811,\,0.146,\,0.231\}$ and
	$e_2(\frac12,\frac12)=\{0.063,\,0.104,\,0.045,\,0.081\}$, respectively.	
	\label{fig:Fig1}}
	\end{center}
	\end{figure}			        			       	
\end{exmp}

In Examples~\ref{Ex:2} and \ref{Ex:3}, we compared our approach with the recently 
published methods 
of Huang \textit{et al.} \cite{HSL08}  and Lu \cite{Lu11}.
The idea of the first method is the following. Given the rational curve \eqref{E:Rn},
use degree elevation to obtain 
\[
	R_n(t)=\frac{\displaystyle  \sum_{i=0}^{n+h} \omega^{<h>}_ir^{<h>}_i
	\,B^{n+h}_i(t)}{ \displaystyle \sum_{i=0}^{n+h} \omega^{<h>}_i\,B^{n+h}_i(t)}
	\qquad (h=0,1,\ldots),
\]
and define the sequence of polynomial curves
\[
	U_{n+h}(t):=\sum_{i=0}^{n+h}r^{<h>}_i\,B^{n+h}_i(t)
	\qquad (h=0,1,\ldots).
\]
Then $\frac{\der^j}{\dt^j}U_{n+h}$ converges uniformly to  $\frac{\der^j}{\dt^j}R_{n}$ as 
$h\to\infty$, 
for any integer $j\ge0$. The weakness of this approach is that the convergence may be
rather slow. Also, notice the increasing degree of the approximating curves.

In the iterative method of Lu \cite{Lu11}, the sequence of B\'ezier curves $\{V^{h}_n\}$
is constructed, 
\[
	V^{h}_n(t)=\sum_{i=0}^{n}v^{h}_i\,B^{n}_i(t)\qquad(h=0,1,\ldots),
\]
where    
\[
	v^{0}_i=R_n(t_i)\qquad(i=0,1,\ldots,n),
\]
with $0=t_0<t_1<\ldots<t_n=1$, and
\[
	v^{h+1}_i=v^{h}_i+\lambda \left(v^{0}_i-V^{h}_n(t_i)\right) 
	\qquad(i=0,1,\ldots,n;\;h=0,1,\ldots),
\]
$\lambda$ being a parameter. It is shown that
\[
	\lim_{h\to\infty}V^{h}_n(t_i)=R_n(t_i)\qquad(i=0,1,\ldots,n).
\]
Also this process may be slow, even for carefully chosen factor $\lambda$
(cf.  \cite{Lu10}). Another drawback of both methods is that only the simplest constraints 
(corresponding to $k=l=1$) are accepted.

In the next two examples  we let $\alpha=\beta=0$.

\begin{exmp}\label{Ex:2}
Our second test rational curve, which has the control points 
 $(14, 1)$,  $(34, 25)$,  $(40, 38)$,  $(-12, 24)$,  $(5, 21)$,  
        $(26, 7)$,  $(18, 41)$,  $(-13, 34)$,  $(14, 1)$
and the associated weights $1,\,  3,\,  3,\,  4,\,  1,\,  7,\,  5,$  $3,\,  1$,
is shown in left part of Figure~\ref{fig:Fig2}  (solid and red).  
We compared our approach with the  methods 
of Huang \textit{et al.} and Lu. 
We produced polynomial approximation of degree $m=10$ with end-point 
interpolation constraints (in our method we put $k=l=1$). 
 The results are shown in the left part of Figure~\ref{fig:Fig2}. Table~\ref{tab:Ex23} 
 below lists the errors
 of approximation in each case (with the number \textit{iter} of iterations performed in 
 the method 
 of  \cite{Lu11}). We see that the method of Huang \textit{et al.}, which is very simple and
 fast, gives  much worse results than two other methods. 
 The convergence of the method of Lu is slow: the result obtained after 100
iterations is 4 times less adequate than ours; also, the comparison of the execution times 
 shows the advantage of our algorithm. 
 \begin{figure}[h!]  
	\begin{center}
        \hspace*{0.5cm} 	
	\includegraphics[scale=0.3,angle=0]{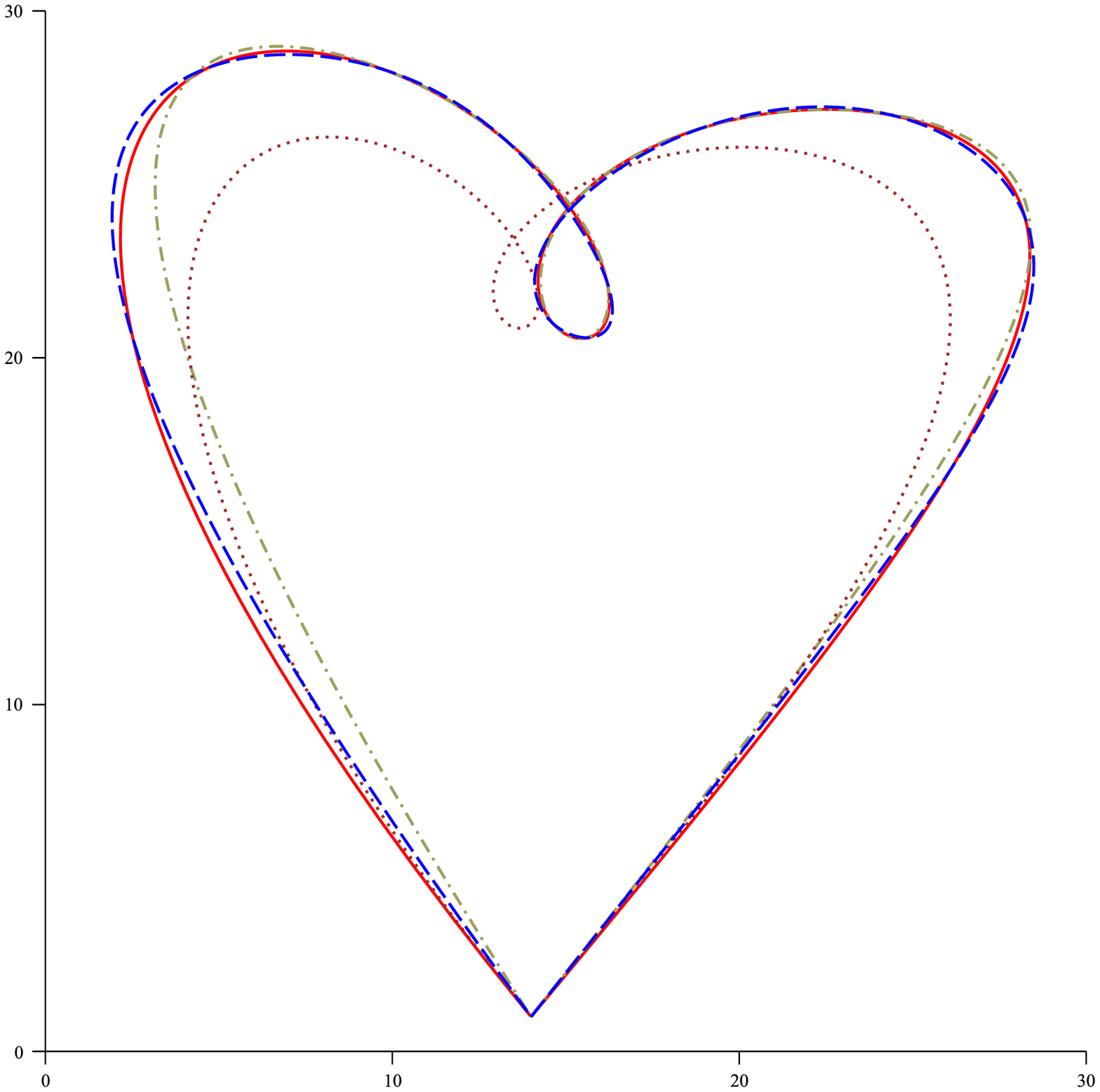}  
	\includegraphics[scale=0.3,angle=0]{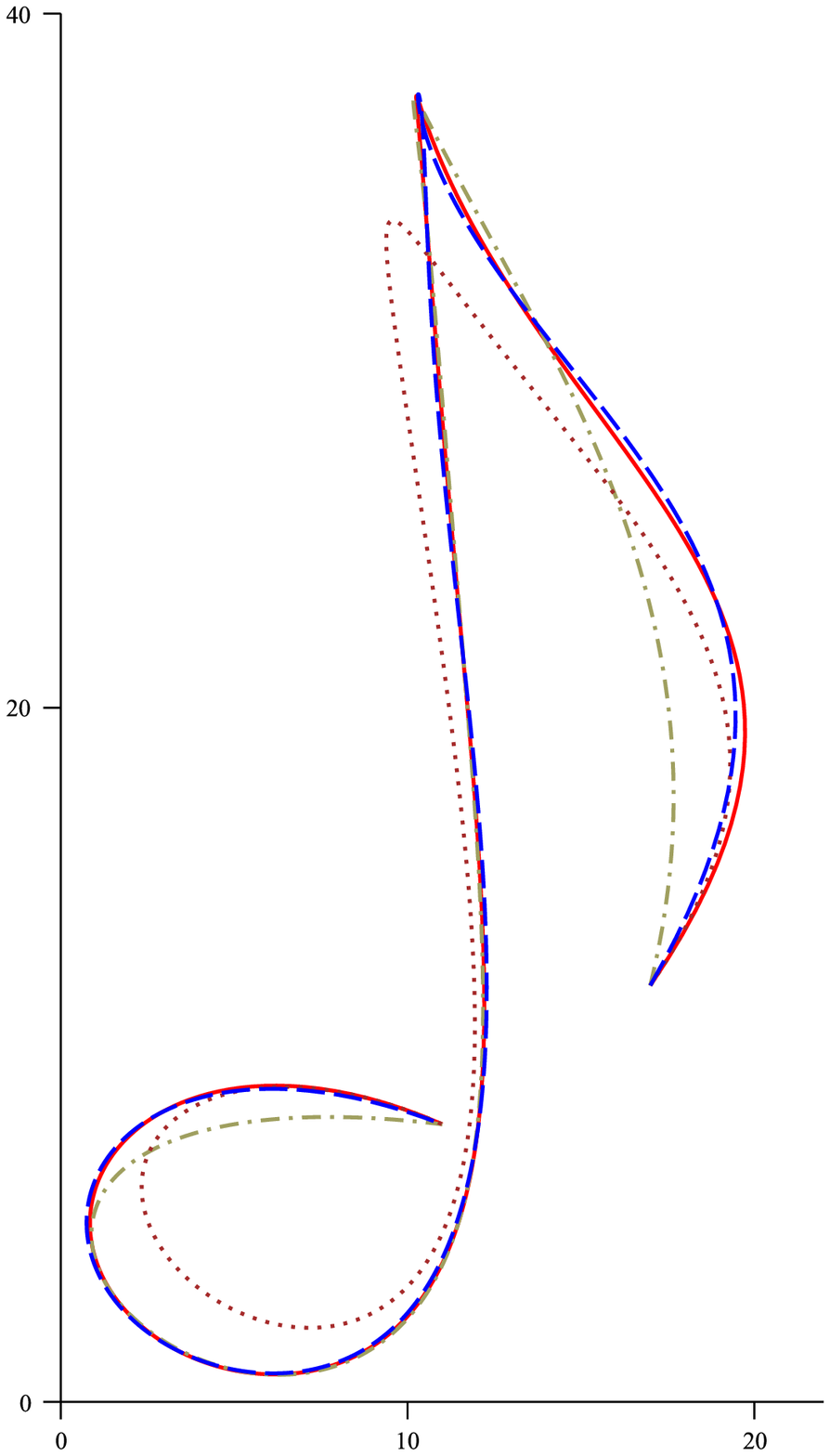}  
        \caption{Polynomial approximation of degree 10
        to the rational B\'ezier curve 
        (solid and red) of degree 8 (\textit{left}) and	
	9 (\textit{right}), with end-point interpolation ($k=l=1$). 
	Comparison of the present method (dashed and  blue), 
        Huang \textit{et al.} method (dotted and brown)
        and Lu method ($100th$ iterate; dotted-dashed and khaki). 		
	\label{fig:Fig2}}
	\end{center}
	\end{figure}
\begin{table}[h!]  
\caption{Errors in the polynomial approximations to the rational B\'ezier
	curve, with end-point interpolation
	\label{tab:Ex23}}
 \small
\renewcommand{\arraystretch}{1.3} 
\[
\begin{array}{lrllr|rllr}
   &      \multicolumn{4}{c}{\sf Example~\ref{Ex:2}} &  \multicolumn{4}{|c}{\sf Example~\ref{Ex:3}}\\
   & iter & e_\infty &\!\!e_2(0,0)& time[s]& iter & e_\infty &\!\! e_2(0,0)& time[s]\\[0.5ex]\hline 
\text{\sf Present method} & \-- & 0.664 & 0.167 & 0.15
                          & \-- & 0.398 & 0.106 & 0.15 \\[1ex]\hline 
\text{\sf Huang \textsl{et al.} \cite{HSL08}}& \-- &  9.41 & 3.98 &<0.01 
                                             & \-- &  5.78 & 3.03 &<0.01 \\[1ex]\hline 
& 25 &  3.43 & 0.850 &0.18& 25 &  2.82 & 0.653 &0.20 \\ 
\text{\sf Lu \cite{Lu11}} &     50 & 2.99 & 0.729&0.35 &     50 & 2.55 & 0.584 &0.37  \\ 
&      75 & 2.79 & 0.671&0.50&      75 & 2.42 & 0.550 &0.54\\ 
&     100 & 2.67 & 0.636 & 0.67&     100 & 2.34 & 0.529 & 0.73 
\end{array}
\]
\end{table}	
\end{exmp}

\begin{exmp}\label{Ex:3}
   $(17, 12)$, $(32, 34)$, $(-23, 24)$, $(33, 62)$, $(-23, 15)$, 
             $(25, 3)$, $(30, -2)$, $(-5, -8)$, $(-5, 15)$, $(11,8)$, 
and the associated weights $1,\, 2,\, 3,\, 6$, $4,\, 5,\, 3,\, 4,\, 2,\, 1$,
is shown in the right part of  Figure~\ref{fig:Fig2}  (solid and red).  
We produced polynomial approximation of degree $m=10$ with end-point 
interpolation constraints. 
 The results are shown in the right part of Figure~\ref{fig:Fig2} and in Table~\ref{tab:Ex23}.
Again, we see  that the methods of Huang \textit{et al.} and Lu give much less
adequate results  than our method.
\end{exmp}

\section*{Conclusions}
                    \label{sec:Concl}

We present a method to solve the 
constrained least squares  approximation
 of the rational B\'ezier curve  by the B\'ezier curve. 
Important tools used are efficient evaluation of  the B\'ezier coefficients
of the constrained  dual Bernstein basis polynomials associated with the Jacobi scalar product
and numerical computation of some integrals involving rational functions.
The new algorithm is particularly attractive when
it is combined with the subdivision process.


\begin{thebibliography}{99}
\itemsep0.5pt
                
\bibitem{Ahn03} Ahn, Y.J.: Using Jacobi polynomials for degree reduction of B\'ezier
                curves with $C^k$-constraints. Comput. Aided Geom. Design \textbf{20},
                423--434 (2003)
                
\bibitem{CW10} Cai, H.J.,  Wang, G.J.: 
               Constrained approximation of rational Bézier curves based on a matrix               
               expression of its end points continuity condition.
               Comput. Aided Design \textbf{42},  495–504 (2010)                       
		
\bibitem{DB08} G.~Dahlquist, A.~Bj\"orck, {Numerical Methods in Scientific
		Computing}, Vol. I, SIAM, 2008.

\bibitem{Far96} G.\,E. Farin,
                {Curves and Surfaces for Computer-Aided
                Geometric Design. A Practical Guide}, third ed.,
                Academic Press,  Boston, 1996.

\bibitem{Flo06}	M.S. Floater, High order approximation of rational curves by polynomial
		curves, Comput. Aided  Design 23 (2006) 621--628.

\bibitem{Gen72} W.M. Gentleman, Implementing Clenshaw-Curtis quadrature. I and II, 
		Comm. ACM 15 (1972) 337--346.
		
\bibitem{HSL08} Y. Huang, H. Su, H. Lin, A simple method for approximating rational
		B\'ezier curve using B\'ezier curves, Comput. Aided Geom. Design
		25 (2008) 697--699.
		 
\bibitem{Kel07} P. Keller, A method for indefinite integration of oscillatory and singular
		functions, Numer. Algor. 46 (2007) 219--251.

\bibitem{KLS10} R.~Koekoek, P. Lesky, R.\,F.~Swarttouw,
               Hypergeometric Orthogonal Polynomials and Their $q$-Analogues,
               Springer, Berlin, 2010.
                    

\bibitem{LP98}  B.-G. Lee, Y. Park, Approximate conversion of rational B\'ezier curves,
                J. KSIAM 2 (1998) 88--93.

\bibitem{Lew92}  S. Lewanowicz, Quick construction of recurrence relations for the Jacobi
	      coefficients,  J. Comput. Appl. Math. 43 (1992) 355-372.		
		
\bibitem{LW11a}  S. Lewanowicz,  P. Wo\'zny,   Multi-degree reduction of tensor product
		B\'ezier surfaces with general constraints,
		{Appl. Math. Comput.} 217 (2011), 4596-4611.
       
\bibitem{LW11b}   S. Lewanowicz,  P. Wo\'zny,   
                B\'ezier representation of the constrained dual Bernstein polynomials,
               {Applied Mathematics and Computation}
	       {218} (2011),  4580--4586.	           
                
\bibitem{LWK12} S. Lewanowicz, P. Wo\'zny,   P. Keller,
                Polynomial  approximation of rational B\'ezier curves with constraints,
                Numer. Algor. 59 (2012) 607--622.	 
		
\bibitem{LW00}   L.G. Liu, G.J. Wang, Recursive formulae for Hermite polynomial approximation
		to rational B\'ezier curves, in: R. Martin, W.P. Wang (Eds.), 
		Proc. Geom. Modeling and Processing 2000: Theory and Applications, 
		IEEE Computer Soc., Los Alamitos, 2000, pp. 190--197.
		
\bibitem{Lu10}	L. Lu, Weighted progressive iteration approximation of rational curves,
                {Comput. Aided Geom. Design}  27 (2010) 127--137.
                
\bibitem{Lu11}	L. Lu, Sample-based polynomial approximation of rational B\'ezier
		curves, J. Comput. Appl. Math. 235 (2011) 1557--1563. 
		
\bibitem{SK91}  T.W. Sederberg, M. Kakimoto, Approximating rational curves using 
		polynomial curves,
		in: G. Farin (Ed.), NURBS for Curve and Surface Design, SIAM, 
		Philadelphia, 1991,
		pp. 144--158.
		
\bibitem{WSC97} G.J. Wang, T.W. Sederberg, F.L. Chen, On the convergence of polynomial
		approximation of rational functions, J. Approx. Theory 89 (1997) 267--288.		
		
\bibitem{WL09}  P. Wo\'zny,  S. Lewanowicz,  Multi-degree reduction of B\'ezier curves 
		with constraints, using dual Bernstein basis polynomials, 
		  {Comput. Aided Geom. Design} 26 (2009) 566--579.
		  
\end{thebibliography}
\end{document}